\documentclass[11pt,leqno,a4paper]{article}

\usepackage{amsmath, amssymb, amsthm, amsfonts}
\usepackage{mathrsfs}
\usepackage{color}

\theoremstyle{plain}
\newtheorem{thm}{Theorem}

\newtheorem{lem}{Lemma}

\theoremstyle{remark}
\newtheorem{rem}{Remark}

\numberwithin{equation}{section}

\begin{document}

\title{Multiple-correction  and Faster Approximation}
\author{Xiaodong Cao, Hongmin Xu and Xu You}
\date{}

\maketitle


\footnote[0]{This work is supported by the National Natural Science Foundation of China (Grant No.11171344) and the Natural
Science Foundation of Beijing (Grant No.1112010).}

\begin{abstract}
In this paper, we formulate a new \emph{multiple-correction method}. The goal is to accelerate the rate of convergence. In particular, we construct some sequences to approximate the Euler-Mascheroni and Landau constants, which are faster than the classical approximations in literature.
\end{abstract}

\section{Introduction}
Euler constant was first introduced by Leonhard Euler (1707-1783) in 1734 as the limit of the sequence
\begin{align}
\gamma(n):=\sum_{m=1}^{n}\frac 1m -\ln n.
\end{align}
It is also known as the Euler-Mascheroni constant. There are many famous unsolved problems about the nature of this constant. For example, it is a long-standing open problem if it is a rational number. See e.g. the survey papers or books of Brent and Zimmermann~\cite{BZ}, Dence and Dence~\cite{DD}, Havil~\cite{Hav} and Lagarias~\cite{Lag}.  A good part of its mystery comes from the fact that the known algorithms converging to $\gamma$ are not very fast, at least, when they are compared to similar algorithms for $\pi$ and $e$.

The sequence $\left(\gamma(n)\right)_{n\in \mathbb{N}}$ converges very slowly toward $\gamma$, like $(2n)^{-1}$. To evaluate it more
accurately, we need to accelerate the convergence. This can be done using the Euler-Maclaurin summation formula, Stieltjes approach
, exponential integral methods, Bessel function method, etc.
See e.g. Gourdon and Sebah~\cite{GS-1}.

Up to now, many authors are preoccupied to improve its rate of convergence. See e.g. Chen and Mortici~\cite{CM}, DeTemple~\cite{De}, Gavrea and Ivan~\cite{GI}, Lu~\cite{Lu,Lu1}, Mortici~\cite{Mor1}, Mortici and Chen~\cite{MC} and references therein. We list some main results as follows: as $n\rightarrow\infty$,
\begin{align}
&\sum_{m=1}^{n}\frac 1m -\ln \left(n+\frac 12\right)=\gamma+O(n^{-2}),
\quad \mbox{(DeTemple~\cite{De}, 1993)},\\
&\sum_{m=1}^{n}\frac 1m -\ln
\frac{n^3+\frac 32n^2+\frac{227}{240}+\frac{107}{480}}
{n^2+n+\frac{97}{240}}=\gamma+O(n^{-6}), \quad \mbox{(Mortici~\cite{Mor1}, 2010)},\\
&\sum_{m=1}^{n}\frac 1m -\ln\rho(n)=\gamma+O(n^{-5}), \quad \mbox{(Chen and Mortici~\cite{CM}, 2012)},
\end{align}
where $\rho(n)=1+\frac {1}{2n}+\frac{1}{24n^2}-\frac{1}{48n^3}+\frac{23}{5760n^4}$.
Recently, Mortici and Chen~\cite{MC} provided a very interesting sequence
\begin{align*}
\nu(n)=&\sum_{m=1}^{n}\frac 1m -\frac 12\ln\left(n^2+n+\frac 13\right)\\
&-\left(\frac {-\frac{1}{180}}{\left(n^2+n+\frac 13\right)^2}
+\frac {\frac{8}{2835}}{\left(n^2+n+\frac 13\right)^3}
+\frac {\frac{5}{1512}}{\left(n^2+n+\frac 13\right)^4}
+\frac {\frac{592}{93555}}{\left(n^2+n+\frac 13\right)^5}
\right),
\end{align*}
and proved
\begin{align}
\lim_{n\rightarrow \infty}n^{12}\left(\nu(n)-\gamma\right)=-\frac{796801}{43783740}.
\end{align}
Hence, the rate of convergence of the sequence $\left(\nu(n)\right)_{n\in \mathbb{N}}$ is $n^{-12}$.

Let $R_1(n)=\frac{a_1}{n}$ and for $k\ge 2$
\begin{align}
R_k(n):=\frac{a_1}{n+\frac{a_2 n}{n+\frac{a_3n}{n+\frac{a_4n}{\frac{\ddots}{n+a_k}}}}
},
\end{align}
where
$(a_1,a_2,a_4,a_6,a_8,a_{10},a_{12})=\left(
\frac{1}{2},\frac{1}{6},\frac{3}{5},
\frac{79}{126},\frac{7230}{6241},
\frac{4146631}{3833346},
\frac{306232774533}{179081182865}\right)$,
$a_{2k+1}=-a_{2k}$ for $1\le k\le 6$, and
\begin{align}
r_k(n):=\sum_{m=1}^{n}\frac 1m -\ln n-R_k(n).
\end{align}
Lu~\cite{Lu} introduced a continued fraction method to investigate
this problem, and showed
\begin{align}
\frac{1}{120(n+1)^{4}}< r_{3}(n)-\gamma<\frac{1}{120(n-1)^{4}}.
\label{Lu-inequalities}
\end{align}
In fact, Lu~\cite{Lu} determined the constants $a_1$ to $a_4$. Xu and You~\cite{XY} continued Lu's work
to find $a_5,\cdots, a_{13}$ with the help of \emph{Mathematica} software, and obtained
\begin{align}
\lim_{n\rightarrow\infty}n^{k+1}\left(r_k(n)-\gamma\right)=C_k',
\label{XY-results}
\end{align}
where
$(C_1',\cdots,C_{13}')=\left(-\frac{1}{12},-\frac{1}{72},\frac{1}{120},
\frac{1}{200},-\frac{79}{25200},-\frac{6241}{3175200},\frac{241}{105840},
\frac{58081}{22018248},-\frac{262445}{91974960},\right.$

$\left. -\frac{2755095121}{892586949408},\frac{20169451}{3821257440},
\frac{406806753641401}{45071152103463200},
-\frac{71521421431}{5152068292800}\right)$. Moreover, they improved \eqref{Lu-inequalities} to
\begin{align}
&C_{10}'\frac{1}{(n+1)^{11}}<\gamma -r_{10}(n)<C_{10}'\frac{1}{n^{11}},\\
&C_{11}'\frac{1}{(n+1)^{12}}< r_{11}(n)-\gamma<C_{11}'\frac{1}{n^{12}}.
\end{align}
However, it seems difficult for us to find more constants $a_k$. One of the main
reasons is due to the recursive algorithm. The other
reason is that the parameter $a_j$ appears many times in the coefficients of polynomials $P_l(x)$ and $Q_m(x)$, and this causes that
expanding function $\frac{P_l(x)}{Q_m(x)}$ as power series in the terms of $1/x$ needs a huge of computations. To overcome this difficulty, the purpose of this paper is to formulate a new \emph{multiple-correction method} to accelerate the convergence. In addition, we will use this method to study the sharp bounds for the constants of Landau.

The Landau's constants are defined for all integers $n\ge 0$ by \begin{align}
G(n)=\sum_{k=0}^{n}\frac{1}{16^k}\binom {2k}{k}^2.
\label{Landau-C-def}
\end{align}
The constants $G(n)$ are important in complex analysis.
In 1913, Landau~\cite{La} proved that if $f(z)=\sum_{k=0}^{\infty}a_kz^k$ is an analytic function in the unit disc which satisfies
$|f(z)|<1$ for $|z|<1$, then $\left|\sum_{k=0}^{n}a_k\right|\le G(n)$, and that this bound is optimal. Landau~\cite{La} showed that
\begin{align}
G(n)\sim \frac{1}{\pi}\ln n, (n\rightarrow \infty).
\label{Landau-C-1}
\end{align}
In 1930, Watson~\cite{Wat} obtained the following more precise asymptotic
formula
\begin{align}
G(n)\sim \frac{1}{\pi}\ln (n+1)+c_0-
\frac{1}{4\pi (n+1)}+O\left(\frac{1}{n^2}\right),
(n\rightarrow \infty),\label{Landau-C-2}
\end{align}
where
\begin{align}
c_0=\frac{1}{\pi}(\gamma+4\ln 2)=1.0662758532089143543\ldots.
\label{c0-definition}
\end{align}
The work of Watson opened up a novel insight into the asymptotic behavior of the Landau sequences $(G(n))_{n\ge0}$. Inspired by formula \eqref{Landau-C-2}, many authors investigated the
upper and lower bounds of $G(n)$. Some of the main results are listed as follows:
\begin{eqnarray}
~~&&\frac{1}{\pi}\ln (n+1)+1\le G(n)
<\frac{1}{\pi}\ln (n+1)+c_0\quad (n\ge 0),\quad\mbox{(Brutman~\cite{Br},1982)},\\
~~&&\frac{1}{\pi}\ln \left(n+\frac 34\right)+c_0<G(n)
\le\frac{1}{\pi}\ln\left(n+\frac 34\right)+1.0976\quad (n\ge 0),\quad\mbox{(Falaleev~\cite{Fal}, 1991)},\label{Landau-Falaleev}\\
~~&&\frac{1}{\pi}\ln \left(n+\frac 34\right)+c_0<G(n)
<\frac{1}{\pi}\ln\left(n+\frac 34+\frac{11}{192n}\right)+c_0\quad (n\ge 1),\quad\mbox{(Mortici~\cite{Mor4}. 2011)}\label{Landau-Mortici}
\end{eqnarray}
Very recently, Chen~\cite{Ch2} presented the following better
approximation to $G(n)$: as $n\rightarrow\infty$,
\begin{align}
G(n)
=&c_0+\frac{1}{\pi}\ln\left(n+\frac 34+\frac{11}{192(n+\frac 34)}-\frac{2009}{184320(n+\frac 34)^3}
+\frac{2599153}{371589(n+\frac 34)^5}\right)\\
&+O\left(\frac{1}{(n+\frac 34)^8}\right),\nonumber
\end{align}
and the better upper bound:
\begin{align}
G(n)
<c_0+\frac{1}{\pi}\ln\left(n+\frac 34+\frac{11}{192(n+\frac 34)}\right), (n\ge 0).\label{chen-1}
\end{align}

Another direction for developing the approximation to $G(n)$ was
initiated by Cvijovi\'c and Klinowski~\cite{CK}, who established estimates
for $G(n)$ in terms of the Psi(or Digamma) function $\psi(z):=\frac{\Gamma'(z)}{\Gamma(z)}$:
\begin{align}
&\frac{1}{\pi}\psi\left(n+\frac 54\right)+c_0<G(n)
<\frac{1}{\pi}\psi\left(n+\frac 54\right)+1.0725,\quad(n\ge 0),\\
&\frac{1}{\pi}\psi\left(n+\frac 32\right)+0.9883<G(n)
<\frac{1}{\pi}\psi\left(n+\frac 32\right)+c_0,\quad(n\ge 0).
\end{align}
Since then, many authors have made significant contributions to sharper
inequalities and asymptotic expansions for $G(n)$. See e.g.
Alzer~\cite{Al}, Chen~\cite{Ch1}, Cvijovi\'c and Srivastava~\cite{CS}, Granath~\cite{Gra}, Mortici~\cite{Mor4}, Nemes~\cite{Ne1,Ne2}, Popa~\cite{Po}, Popa and Secelean~\cite{PS}, Zhao~\cite{Zhao}, Gavrea and M. Ivan~\cite{GI}, Chen and Choi~\cite{CC1,CC2}, etc. To the best knowledge of authors, the latest upper bound is due to Chen~\cite{Ch2}, who proved
\begin{align}
G(n)
<c_0+\frac{1}{\pi}\psi\left(n+\frac 54+\frac{1}{64(n+\frac 34)}\right),  (n\ge 0).\label{chen-2}
\end{align}
 Here, the authors would like to thank Alzer, Chen, Choi, DeTemple, Granath, Lu, Mortici, etc., it is their important works that makes the present work  becomes possible.

\noindent {\bf Notation.} Throughout the paper,  the notation
$P_k(x)$(or $Q_k(x)$) as usual denotes a polynomial of degree $k$
in terms of $x$. The notation $\Psi(k;x)$
means a polynomial of degree $k$ in terms of $x$ with all of its non-zero coefficients being positive, which may be different at each occurrence.
Notation $\Phi(k;x)$ denotes a polynomial of degree $k$ in terms of $x$ with the leading coefficient being equal to one, which may be different at different subsection.

\section{Some Lemmas}

The following lemma gives a method for measuring the rate of convergence, for the proof of which, see Mortici~\cite{Mor1,Mor2}.
\begin{lem}
If the sequence $(x_n)_{n\in \mathbb{N}}$ is convergent to zero and there
exists the limit
\begin{align}
\lim_{n\rightarrow +\infty}n^s(x_n-x_{n+1})=l\in [-\infty,+\infty]
\label{LEM1-1}
\end{align}
with $s>1$, then
\begin{align}
\lim_{n\rightarrow +\infty}n^{s-1}x_n=\frac{l}{s-1}.\label{LEM1-2}
\end{align}
\end{lem}

In the study of Landau constants, we need to apply a so-called
Brouncker's continued fraction formula.
%, which plays an important role.
\begin{lem}
For all integer $n\ge 0$, we have
\begin{align}
q(n):=\left(\frac{\Gamma(n+\frac 12)}{\Gamma(n+1)}\right)^2=
\frac{4}{1+4n+\frac{1^2}{2+8n+\frac{3^2}{2+8n
+\frac{5^2}{2+8n+\ddots}}}}.\label{LEM2-1}
\end{align}
\end{lem}

In 1654 Lord William Brouncker found this remarkable fraction formula, when Brouncker and Wallis collaborated on the problem of squaring the circle. Formula (2.3) was not published by Brouncker himself, but first appeared in ~\cite{Wal}. For a general $n$, Formula (2.3) follows from Entry 25 in Chapter 12 in Ramanujan¡¯s notebook ~\cite{Ber}, which gives a more general continued fraction formula for quotients of gamma functions, and which have several proofs published by different
authors.

Writing continued fractions in this way of~\eqref{LEM2-1} takes a lot of space. So instead we use the following shorthand notation
\begin{align}
q(n)=\frac{4}{1+4n+}\frac{1^2}{2+8n+}\frac{3^2}{2+8n+}
\frac{5^2}{2+8n+}\cdots,
\end{align}
and its $k$-th approximation $q_k(n)$ is defined by
\begin{align}
q_1(n)&=\frac{4}{1+4n},\\
q_k(n)&=\frac{4}{1+4n+}\frac{1^2}{2+8n+}\frac{3^2}{2+8n+}\cdots
\frac{(2k-3)^2}{2+8n},\quad (k\ge 2).\label{qk-definition}
\end{align}

In the proof of our inequalities for the constants of Euler-Mascheroni and Landau, we also use the following simple inequality.
\begin{lem} Let $f''(x)$ be a continuous function.  If $f''(x)>0$, then
\begin{align}
\int_{a}^{a+1}f(x) dx > f(a+1/2).\label{LEM3}
\end{align}
\end{lem}
\proof By letting $x_0=a+1/2$ and Taylor's formula, we have
\begin{align*}
\int_{a}^{a+1}f(x) dx=&\int_{a}^{a+1}\left(
f(x_0)+f'(x_0)(x-x_0)+\frac 12 f''(\theta_x)(x-x_0)^2\right)dx\\
>&\int_{a}^{a+1}\left(
f(x_0)+f'(x_0)(x-x_0)\right)dx=f(a+1/2).
\end{align*}
This completes the proof of Lemma 3. Also see Lemma 2 in Xu and You ~\cite{XY}.\qed
\section{Two Examples for Euler-Mascheroni Constant}
In this section, to illustrate quickly and clearly the main ideas of this paper, we consider the simplest case of Euler-Mascheroni constant by using the correction-process again.

\noindent{\bf Example 1.} We choose an initial-correction function $\theta_0(n)$ given by
\begin{align}
\theta_0(n)=
\frac{13 + 30 n}{6 (1 + 6 n + 10 n^2)},\label{example1-1}
\end{align}
and define
\begin{align}
\nu_0(n)=\sum_{m=1}^{n}\frac 1m -\ln n -\theta_0(n).\label{example1-2}
\end{align}
Applying Lemma 1, one can check(See Theorem 1 in ~\cite{XY}
, or \eqref{XY-results})
\begin{align}
&\lim_{n\rightarrow\infty}n^{6}\left(\nu_0(n)-\nu_0(n+1)\right)
=\frac{1}{40},\label{example1-3}\\
&\lim_{n\rightarrow\infty}n^{5}\left(\nu_0(n)-\gamma\right)
=\frac{1}{200}.\label{example1-4}
\end{align}

By using the similar idea of Kummer's acceleration method and inserting the correction function $-\frac{1}{200}\frac{1}{n^5}$ in \eqref{example1-2} again,
one can use Lemma 1 again to show
\begin{align}
\lim_{n\rightarrow\infty}n^{6}
\left(\nu_0(n)-\frac{1}{200}\frac{1}{n^5}-\gamma\right)
=-\frac{773}{126000}.\label{example1-5}
\end{align}
%Indeed, the above result is almost meaningless.
Furthermore, we try to obtain an algorithm with a faster convergent rate by using  $\Phi(5;n)$ instead of $n^5$. To do that, let
\begin{align}
\theta(n)=&\frac{\frac{1}{200}}{\Phi(5;n)}=\frac{\frac{1}{200}}
{n^5+a_4 n^4+a_3 n^3
+a_2 n^2+a_1 n+a_0},\label{example1-6}\\
\nu(n)=&\sum_{m=1}^{n}\frac 1m -\ln n -\theta_0(n)-\theta(n).
\label{example1-7}
\end{align}
First, we use the method of undetermined coefficients to find $a_{j}(0\le j\le 4)$. By using the \emph{Mathematica} software, we
expand the difference $\nu(n)-\nu(n+1)$ into a power series in terms of $n^{-1}$:
\begin{align}
&\nu(n)-\nu(n+1)\label{example1-equation}\\
=&\frac{-773+630 a_4}{21000}\frac{1}{n^7}+
\frac{4033+1050 a_3-3150a_4-1050 a_4^2}{30000}\frac{1}{n^8}\nonumber\\
&+\frac{-37657+4500 a_2-15750 a_3+31500 a_4-9000 a_3 a_4+15750 a_4^2+
4500 a_4^3}{112500}\frac{1}{n^9}\nonumber\\
&+\frac{\varphi_{10}}{500000}\frac{1}{n^{10}}
+\frac{\varphi_{11}}{4125000}
\frac{1}{n^{11}}+\frac{\varphi_{12}}{75000000}\frac{1}
{n^{12}}+O\left(
\frac{1}{n^{13}}\right),\nonumber
\end{align}
where
\begin{align*}
\varphi_{10}=&350191+22500 a_1-90000 a_2+210000 a_3-22500 a_3^2
-315000 a_4-
45000 a_2 a_4\\
&+180000 a_3 a_4-210000 a_4^2+67500a_3 a_4^2-
90000 a_4^3-22500 a_4^4,\\
\varphi_{11}=&-5465923+206250 a_0-928125 a_1+2475000 a_2
-4331250 a_3 -
412500 a_2 a_3\\
&+928125 a_3^2+5197500 a_4-412500 a_1 a_4+
1856250 a_2 a_4
-4950000 a_3 a_4\\
&+618750 a_3^2 a_4+4331250 a_4^2+
618750 a_2 a_4^2-2784375 a_3 a_4^2+2475000 a_4^3\\
&-825000 a_3 a_4^3+
928125 a_4^4+206250 a_4^5,\\
\varphi_{12}=&175990871-20625000 a_0+61875000 a_1-123750000a_2-
4125000 a_2^2\\
&+173250000 a_3-8250000 a_1 a_3+41250000 a_2 a_3-61875000 a_3^2+
4125000 a_3^3\\
&-173250000 a_4-8250000 a_0 a_4+41250000 a_1 a_4-
123750000 a_2 a_4+247500000 a_3 a_4\\
&+24750000 a_2 a_3 a_4-
61875000 a_3^2 a_4-173250000 a_4^2+12375000 a_1 a_4^2-
61875000 a_2 a_4^2\\
&+185625000 a_3 a_4^2-24750000 a_3^2 a_4^2-
123750000 a_4^3-16500000 a_2 a_4^3\\
&+82500000 a_3 a_4^3-
61875000 a_4^4+20625000 a_3 a_4^4-20625000 a_4^5-4125000 a_4^6.
\end{align*}

According to Lemma 1, we have five parameters $a_4, a_3, a_2, a_1$ and $a_0$ which produce the fastest convergence of the sequence
from \eqref{example1-equation}
\begin{align*}
\begin{cases}
-773+630 a_4=0,\\
4033+1050 a_3-3150a_4-1050 a_4^2=0\\
-37657+4500 a_2-15750 a_3+31500 a_4-9000 a_3 a_4+15750 a_4^2+
4500 a_4^3=0\\
\varphi_{10}=0,\\
\varphi_{11}=0,
\end{cases}
\end{align*}
namely if
\begin{align}
\theta(n)=\frac{\frac{1}{200}}{
n^5+\frac{773}{630}n^4+\frac{21361}{15876}n^3
+\frac{1348075}{2000376}n^2-
\frac{91207415}{252047376}n-\frac{178345771979}{1746688315680}
}.\label{example1-8}
\end{align}
Thus, we get
\begin{align}
\nu(n)-\nu(n+1)
=-\frac{10992878936527}{160060165655040}\frac{1}{n^{12}}+O\left(
\frac{1}{n^{13}}\right).\label{example1-9}
\end{align}
We can apply another approach to find $a_4, a_3, a_2, a_1$ and $a_0$ step by step, which is achieved by using  $n^5+a_{4}n^4, n^5+a_{4}n^4+a_{3}n^3, n^5+a_{4}n^4+a_{3}n^3+a_{2}n^2, n^5+a_{4}n^4+a_{3}n^3+a_{2}n^2+a_{1}n$, $ n^5+a_{4}n^4+a_{3}n^3+a_{2}n^2+a_{1}n+a_{0}$ instead of $\Phi(5;n)$ in turn.  For the reader's convenience, here we give an example to explain how
\emph{Mathematica } software generates $\nu(n)-\nu(n+1)$ into power series in the terms of $\frac 1n$. For example, find $a_3$. We manipulate \emph{Mathematica} program
\begin{align*}
&\mbox{Normal[Series[}\left(-\frac{1}{n + 1} + \log[1 + \frac 1n]-\left(\frac{13+30n}{6(1+6 n
+10 n^2)}+\frac{1/200}{n^5+(773/630)n^4+ a_3 n^3}\right)\right.\\
&\left.+\left(\frac{13+30n}{6(1+6 n+10 n^2)}+\frac{1/200}{n^5+(773/630)n^4+ a_3 n^3}\right)/.n\rightarrow n+1\right)/.n\rightarrow 1/x, \{x, 0, 10\}]]/.x\rightarrow 1/n
\end{align*}
to generate
\begin{align}
\nu(n)-\nu(n+1)
=\frac{-\frac{21361}{453600}+\frac{7 a3}{200}}{n^8}+
\frac{\frac{9173092}{31255875}-\frac{3751 a3}{15750}}{n^9}+
O\left(
\frac{1}{n^{10}}\right).\label{example1-equation-1}
\end{align}
By solving the equation $-\frac{21361}{453600}+\frac{7 a3}{200}=0$, we also find $a_3=\frac{21361}{15876}$. In what follows, we always use this approach.

By Lemma 1 again, we obtain finally
\begin{align}
\lim_{n\rightarrow\infty}n^{11}\left(\nu(n)-\gamma\right)
=-\frac{10992878936527}{1760661822205440}.\label{example1-10}
\end{align}
We observe that the above twice-correction improves the rate of convergence from $n^{-6}$ to $n^{-11}$.\qed

\begin{rem}
The main idea of twice-correction is that from $n^5, n^5+a_{4}n^4, n^5+a_{4}n^4+a_{3}n^3+a_{2}n^2, n^5+a_{4}n^4+a_{3}n^3+a_{2}n^2+a_{1}n$ to $ n^5+a_{4}n^4+a_{3}n^3+a_{2}n^2+a_{1}n+a_{0}$, their approximations in turn become better and better.
\end{rem}
\begin{rem}
It should be noted that once we find the exact values of the parameters $a_4$ to $a_0$, it is not very difficult for us to check the formula
\eqref{example1-9} with the help of \emph{Mathematica} software.
\end{rem}

\noindent{\bf Example 2.} We would like to give another example. Now we take the initial-correction function $\eta_0(n)=\frac{6n-1}{12n^2}$(see Theorem 1.1 in Lu~\cite{Lu} or \eqref{XY-results}, which is found by the continued fraction method), and define
\begin{align}
\mu_0(n)=&\sum_{m=1}^{n}\frac 1m -\ln n -\eta_0(n).\label{example2-1}
\end{align}
One may check by using Lemma 1
\begin{align}
\lim_{n\rightarrow\infty}n^{4}\left(\mu_0(n)-\gamma\right)
=\frac{1}{120}.\label{example2-2}
\end{align}
Similarly, we insert a correction function $-\eta(n)$ in
\eqref{example2-1} again, which has the form of
\begin{align}
\eta(n)=\frac{\frac{1}{120}}{n^4+b_3 n^3+b_2 n^2+b_1 n+b_0}.\label{example2-3}
\end{align}
 Applying the same method as Example 1, we can find $(b_3,b_2,b_1,b_0)=(0,\frac{10}{21},0,-\frac{241}{882})$. The details is omitted here. Now we define
\begin{align}
\mu(n)=&\sum_{m=1}^{n}\frac 1m -\ln n -\eta_0(n)-\eta(n),\label{example2-4}
\end{align}
where
\begin{align}\eta_0(n)=\frac{6n-1}{12n^2}\quad \mbox{and}\quad
\eta(n)=\frac{\frac{1}{120}}{n^4+\frac{10}{21}n^2-\frac{241}{882}}.
\label{example2-5}
\end{align}
By using \emph{Mathematica} software and Lemma 1, we can attain
\begin{align}
&\mu(n)-\mu(n+1)
=-\frac{13775}{305613} \frac{1}{n^{11}}+O\left(
\frac{1}{n^{12}}\right),\label{example2-6}\\
&\lim_{n\rightarrow\infty}n^{10}\left(\mu(n)-\gamma\right)
=-\frac{13775}{3056130}.\label{example2-7}
\end{align}

\begin{rem} We observe that the above twice-correction improves the rate of convergence from $n^{-4}$ to $n^{-10}$, which is the desired result. However, it is interesting to note that both $b_3$ and  $b_1$ equal zero. The reason of why
inserting the sub-correction term $b_3n^3$(or $b_1n$)
 does not improves the rate of convergence(i.e. compare $n^4+ b_3n^3$ with $n^4$, or $n^4+ b_3n^3+b_2n^2+b_1n$ with $n^4+ b_3n^3+b_2n^2$)  may be that the function $n^3$(or $n$) \emph{changes} too rapidly when $n$ tends to infinity.
 Fortunately, these losses are made up by the sub-correction
 terms $b_2n^2$ and $b_0$.
%(oscillates)
\end{rem}
More precisely, we will improve \eqref{example2-7}, and prove the following double-sides inequalities.
\begin{thm} Let $\mu(n)$ be defined by \eqref{example2-4}. Then
for all positive integer $n$, we have
\begin{align}
\frac{13775}{3056130}\frac{1}{(n+\frac 34)^{10}}<\gamma-\mu(n)<\frac{13775}{3056130}\frac{1}{(n-\frac 14)^{10}}.\label{example2-thm1}
\end{align}
\end{thm}
\begin{rem}
In fact, Theorem 1 implies that
$\mu(n)$ is a strictly increasing function of $n$.
\end{rem}
\proof
It follows from \eqref{example2-4}
\begin{align}
\mu(n)-\mu(n+1)=-\frac{1}{n + 1}+\ln(1+\frac 1n)+\eta_0(n+1)+\eta(n+1)-\eta_0(n)-\eta(n).\label{example2-8}
\end{align}
We write $D=\frac{13775}{27783}$, and define for $x\ge 1$
\begin{align}
-\omega(x)=-\frac{1}{x + 1}+\ln(1+\frac 1x)+\eta_0(x+1)+\eta(x+1)-\eta_0(x)-\eta(x).\label{example2-9}
\end{align}
By \emph{Mathematica} software, it is not difficult to check
\begin{align}
&-\omega'(x)-\frac{D}{(x+\frac 14)^{12}}\label{example2-10}\\
=&-\frac{\Psi_1(20;x)(x-1)+4032098201877889488940287625}{277830 x^3(1+x)^3(1+4x)^{12}(-241+420 x^2+
882 x^4)^2(1061+4368 x+5712 x^2+3528 x^3+882 x^4)^2}\nonumber\\
<&0,\nonumber
\end{align}
and
\begin{align}
&-\omega'(x)-\frac{D}{(x+\frac 34)^{12}}\label{example2-11}\\
=&\frac{\Psi_2(20;x)(x-1)+51726219719747325679290363431}{277830 x^3(1+x)^3(3+4 x)^{12}(-241+420 x^2+
882 x^4)^2(1061+4368 x+5712 x^2+3528 x^3+882 x^4)^2}\nonumber\\
>&0.\nonumber
\end{align}
Note that $\omega(\infty)=0$. From \eqref{example2-10} and Lemma 3, one has
\begin{align}
\omega(n)=&\int_{n}^{\infty}-\omega'(x)dx<D\int_{n}^{\infty}
\frac{dx}{(x+\frac 14)^{12}}\label{example2-12}\\
=&\frac{D}{11}\frac{1}{(n+\frac 14)^{11}}
<\frac{D}{11}\int_{n-\frac 14}^{n+\frac 34}\frac{dx}{x^{11}}.\nonumber
\end{align}
Note that $\mu(\infty)=\gamma$. Combining \eqref{example2-8},\eqref{example2-9} and \eqref{example2-12}, we have
\begin{align}
\gamma-\mu(n)=&\sum_{m=n}^{\infty}\left(\mu(m+1)-\mu(m)\right)
<\frac{D}{11}\sum_{m=n}^{\infty}\int_{m-\frac 14}^{m+\frac 34}\frac{dx}{x^{11}}\label{example2-13}\\
=&\frac{D}{11}\int_{n-\frac 14}^{\infty}\frac{dx}{x^{11}}
=\frac{D}{110}\frac{1}{(n-\frac 14)^{10}}.\nonumber
\end{align}
This finishes the proof of right-hand inequality in \eqref{example2-thm1}. Similarly, it follows from \eqref{example2-11}
\begin{align}
\omega(n)=&\int_{n}^{\infty}-\omega'(x)dx>D\int_{n}^{\infty}
\frac{dx}{(x+\frac 34)^{12}}\label{example2-14}\\
=&\frac{D}{11}\frac{1}{(n+\frac 34)^{11}}
>\frac{D}{11}\int_{n+\frac 34}^{n+\frac 74}\frac{dx}{x^{11}}.\nonumber
\end{align}
Finally, by \eqref{example2-8},\eqref{example2-9} and \eqref{example2-14}, one has
\begin{align}
\gamma-\mu(n)=&\sum_{m=n}^{\infty}\left(\mu(m+1)-\mu(m)\right)
>\frac{D}{11}\sum_{m=n}^{\infty}\int_{m+\frac 34}^{m+\frac 74}\frac{dx}{x^{11}}\label{example2-15}\\
=&\frac{D}{11}\int_{n+\frac 34}^{\infty}\frac{dx}{x^{11}}
=\frac{D}{110}\frac{1}{(n+\frac 34)^{10}}.\nonumber
\end{align}
This completes the proof of Theorem 1.\qed
\section{The multiple-correction method}
Based on the work of Section 2, we will formulate a new \emph{multiple-correction method} to study faster approximation problem for the constants of Euler-Mascheroni and Landau.

 Let $(v(n))_{n\ge 1}$ be a sequence to be approximated. Throughout the paper, we always assume that the following three
conditions hold.

\noindent {Condition (i).} The initial-correction function $\eta_0(n)$ satisfies
\begin{align*}
&\lim_{n\rightarrow\infty}\left(v(n)-\eta_0(n)\right)=0,\\
&\lim_{n\rightarrow\infty}n^{l_0}\left(v(n)-v(n+1)-\eta_0(n)
+\eta_0(n+1)
\right)=C_0\neq 0,
\end{align*}
with some a positive integer $l_0\ge 2$.

\noindent {Condition (ii).} The $k$-th correction function $\eta_k(n)$ has the form of $-\frac{C_{k-1}}{\Phi_k(l_{k-1};n)}$, where
\begin{align*}
\lim_{n\rightarrow\infty}n^{l_{k-1}}\left(v(n)-v(n+1)-
\sum_{j=0}^{k-1}\left(\eta_j(n)-\eta_j(n+1)\right)
\right)=C_{k-1}\neq 0,
\end{align*}

\noindent {Condition (iii).} The difference  $\left(v(1/x)-v(1/x+1)-\eta_0(1/x)
+\eta_0(1/x+1)
\right)$ is an analytic function in a neighborhood of point $x=0$.
%v(x)\in C^{\infty}[1,+\infty]$.

\subsection{Euler-Mascheroni Constant}
\noindent{\bf (Step 1) The initial-correction.}
We choose $\eta_0(n)=0$, and let
\begin{align}
\nu_0(n)=\sum_{m=1}^{n}\frac 1m -\ln n -\eta_0(n)=\sum_{m=1}^{n}\frac 1m -\ln n.\label{example3-1}
\end{align}
By lemma 1, it is not difficult to prove that
\begin{align}
&\lim_{n\rightarrow \infty}n^2\left(\nu_0(n)-\nu_0(n+1)\right)
=\frac 12,\label{example3-2}\\
&\lim_{n\rightarrow \infty}n\left(\nu_0(n)-\gamma\right)=\frac{1}{2}=:C_0.
\label{example3-3}
\end{align}

\noindent{\bf (Step 2) The first-correction.} We let
\begin{align}
\eta_1(n)=\frac{C_0}{\Phi_1(1;n)}
=\frac{\frac 12}{n+b_{(1,0)}},\label{example3-4}
\end{align}
and define
\begin{align}
\nu_1(n)=\sum_{m=1}^{n}\frac 1m -\ln n -\eta_0(n)-\eta_1(n).\label{example3-5}
\end{align}
By the same method as \eqref{example1-equation-1}, we find $b_{(1,0)}=\frac 16$.
Applying Lemma 1 again, one has
\begin{align}
&\lim_{n\rightarrow \infty}n^4\left(\nu_1(n)-\nu_1(n+1)\right)=-\frac{1}{24},
\label{example3-6}\\
&\lim_{n\rightarrow \infty}n^{3}\left(\nu_1(n)-\gamma\right)=-\frac{1}{72}:=C_1.
\label{example3-7}
\end{align}

\noindent{\bf (Step 3) The second-correction.}  Similarly, we set
the second-correction function in the form of $\eta_2(n)=\frac{C_1}{\Phi_2(3;n)}$, and define
\begin{align}
\nu_2(n)=\sum_{m=1}^{n}\frac 1m -\ln n -\eta_0(n)-\eta_1(n)-\eta_{2}(n).\label{example3-8}
\end{align}
By using similar approach of \eqref{example1-equation-1}, we can find
\begin{align}
\eta_2(n)=
\frac{-\frac{1}{72}}{n^3+\frac{23}{30}n^2+\frac{14}{25}n
+\frac{1333}{10500}}=\frac{-\frac{1}{72}}
{(n+\frac{23}{90})^3+\frac{983}{2700}(n+\frac{23}{90})
+\frac{2197}{127575}}.\label{example3-9}
\end{align}
By Lemma 1, one can obtain
\begin{align}
&\lim_{n\rightarrow \infty}n^8\left(\nu_2(n)-\nu_2(n+1)\right)=-\frac{7061}{5400},
\label{example3-10}\\
&\lim_{n\rightarrow \infty}n^{7}\left(\nu_2(n)-\gamma\right)=
-\frac{7061}{3780000}:=C_2.\label{example3-11}
\end{align}

\noindent{\bf (Step 4) The third-correction.} We set
$\eta_3(n)=\frac{C_2}{\Phi_3(7;n)}$, and define
\begin{align}
\nu_3(n)=\sum_{m=1}^{n}\frac 1m -\ln n -\eta_0(n)-\eta_1(n)-\eta_{2}(n)-\eta_{3}(n).\label{example3-12}
\end{align}
By using \emph{Mathematica} software, we can find
\begin{align}
\Phi_3(7;n)=&n^7+\frac{126901}{70610}n^6+
\frac{302657774122}{78525910575}n^5+
\frac{203050524517143511}{52278737145178500}n^4\label{example3-13}\\
&+\frac{4586975399311716291806}{2713180197918474605475}n^3+
\frac{139644786102811402696525439}{298861139889036647352440010}n^2
\nonumber\\
&+\frac{1335669056713380727335512329403306}{
243734857761374337083369364175455}n\nonumber\\
&+
\frac{1768275723433572920281319954725767947341}{
1032607098391838516487402648265732653000}.\nonumber
\end{align}
Now by Lemma 1 again, we obtain
\begin{thm} Let $\nu_3(n)$ be defined by \eqref{example3-12}. Then
we have
\begin{align}
\lim_{n\rightarrow \infty}n^{16}\left(\nu_3(n)-\nu_3(n+1)\right)=
15C_3
,\label{example3-14}
\end{align}
and
\begin{align}
\lim_{n\rightarrow \infty}n^{15}\left(\nu_3(n)-\gamma\right)=
-\frac{6044981017774921659252823535814990412377703}{
102460439337930176798462527774167322493925000}:=C_3.
\label{example3-15}
\end{align}
\end{thm}
\begin{rem} It could be imagined that if we apply the correction-process many times, then, we can obtain $k$th-correction sequence
\begin{align}
\nu_k(n)=\sum_{m=1}^{n}\frac 1m -\ln n -\eta_0(n)-\sum_{j=1}^{k}
\frac{C_{j-1}}{\Phi_j(l_{j-1};n)}\label{example3-16}
\end{align}
with the rate of convergence $\frac{1}{n^{l_k}}$,
here $l_k\ge 2l_{k-1}+1$, i.e.
\begin{align}
&\lim_{n\rightarrow \infty}n^{l_k}\left(\nu_k(n)-\gamma\right)={C_k}\neq 0,\\
&\gamma=\sum_{m=1}^{n}\frac 1m -\ln n -\eta_0(n)-\sum_{j=1}^{k}
\frac{C_{j-1}}{\Phi_j(l_{j-1};n)}+O\left(\frac{1}{n^{l_k}}\right).
\end{align}
\end{rem}
\begin{rem}
For comparison, the result $\nu_1(n)$  in Theorem 2 is the same as $r_2(n)$ in~\eqref{XY-results}, and $\lim_{n\rightarrow\infty}\frac{\nu_2(n)-\gamma}{r_6(n)-\gamma}
=\frac{C_2}{C_6'}=0.950367\cdots<1$.
Theoretically at least, for a large $n$ the above formula may reduce or eliminate
numerically computations compared with Euler-Maclaurin summation formula. For example, if we take $n=2^{15}=32768$ in Theorem 2, then $ -1.09418\cdot 10^{-69}<\nu_3(n)-\gamma<0$.
\end{rem}
%\begin{rem}
%Theorem 2 predicts that it may be possible for us to find a
% BBP-type series of Euler-Mascheroni constant in the future.
%\end{rem}
\begin{rem}  We can take different initial-correction function to find some other simple faster approximations. For example, we choose the initial-correction function $\eta_0(n)=-\ln n+\frac 12\ln\left(n^2+n+\frac {1}{3}\right)$, see, Chen and Li~\cite{CL}. By Lemma 1, it is not very difficult for us to check that
 \begin{align}
&\nu_0(n)=\sum_{m=1}^{n}\frac 1m -\ln n -\eta_0(n)=\sum_{m=1}^{n}\frac 1m -\frac 12\ln\left(n^2+n+\frac 13\right),\\
&\lim_{n\rightarrow \infty}n^4\left(\nu_0(n)-\gamma\right)=-\frac{1}{180}=:C_0.
\end{align}
Let
\begin{align}
\eta_1(n)=\frac{C_0}{\Phi_1(4;n)}
\quad \mbox{and}\quad
\eta_2(n)=\frac{C_1}{\Phi_2(10;n)},
\end{align}
where $C_1=\frac{457528}{123773265}$,
\begin{align}
\eta_1(n)=&(n+\frac 12)^4+\frac{85}{126}(n+\frac 12)^2-\frac{18287}{63504}
,\\
\eta_2(n)=&(n+\frac 12)^{10}+\frac{28038237821}{5995446912}(n+\frac 12)^8+\frac{11612938185382451401}{35945383674610335744}(n+\frac 12)^6
\\
&+\frac{163544744039006129564874642269}{8307589279805355451415003136}
(n+\frac 12)^4\nonumber\\
&-\frac{2762081970439756978947606523226093660107}{
15021373006058621011789214048600457216}(n+\frac 12)^2\nonumber\\
&+\frac{2771007475606973680958970352585491511233080640189551}{
1350897666047614749541649384438829437061829754880}.\nonumber
\end{align}
Now define
\begin{align}
\nu_2(n)=\sum_{m=1}^{n}\frac 1m -\ln n -\eta_0(n)-\eta_1(n)-\eta_{2}(n).
\end{align}
By using Lemma 1, one may check
\begin{align}
\lim_{n\rightarrow \infty}n^{22}\left(\nu_2(n)-\gamma\right)
={C_2},
\end{align}
where
\begin{align}
C_2=\frac{186484155415412379058939871115843737323085706271910
2654146029}{
18653565767176841210967548892254397636853629986414159462400}.
\end{align}
Some other interesting correction functions can be found in Gourdon and Sebah~\cite{GS}.
\end{rem}

\subsection{Landau Constants}

\noindent{\bf (Step 1) The initial-correction.}
Let $c_0$ be defined by \eqref{c0-definition}. Motivated by inequalities \eqref{Landau-Falaleev} and \eqref{Landau-Mortici},
we choose $\eta_0(n)=\frac{1}{\pi}\ln(n+\frac 34)+c_0$, and define
\begin{align}
u_0(n)=G(n)-\eta_0(n)=G(n)-\frac{1}{\pi}\ln(n+\frac 34)-c_0.
\label{L-u0-definition}
\end{align}
Now we consider the difference $u_0(n)-u_0(n+1)$. It follows
immediately from \eqref{L-u0-definition}
\begin{align}
u_0(n)-u_0(n+1)=G(n)-G(n+1)-\frac{1}{\pi}\ln(n+\frac 34)+
\frac{1}{\pi}\ln(n+\frac 74).\label{L-u0-difference}
\end{align}
First, from the duplication formula (Legendre, 1809)
\begin{align}
2^{2z-1}\Gamma(z)\Gamma(z+\frac 12)=\sqrt\pi\Gamma(2z),
\label{duplication formula}
\end{align}
one can prove
\begin{align}
G(n)-G(n-1)=\frac{(\Gamma(2n+1))^2}{16^2(\Gamma(n+1))^4}
=\left(\frac{(2n)!}{4^n(n!)^2}\right)^2=\frac{1}{\pi}q(n),
\label{G-difference-def}
\end{align}
where $q(n)$ is defined by \eqref{LEM2-1}. Also see p.739 in Granath~\cite{Gra} or p.306 in Chen~\cite{Ch2}.
By \eqref{L-u0-difference} and \eqref{G-difference-def}, one has
\begin{align}
u_0(n)-u_0(n+1)=-\frac{1}{\pi}q(n+1)-\frac{1}{\pi}\ln(n+\frac 34)+
\frac{1}{\pi}\ln(n+\frac 74).\label{L-u0-difference-1}
\end{align}
From Lemma 2 and \eqref{qk-definition}, on one hand, it can be observed that for all positive integer $j$, one has
\begin{align}
q_2(n)<q_4(n)<\cdots<q_{2j}(n)<q(n)
<q_{2j+1}(n)<\cdots<q_3(n)<q_1(n).\label{Landau-1}
\end{align}
On the other hand, we can check by using \emph{Mathematica} software
\begin{align}
q_9(n)-q_8(n)=O\left(\frac{1}{n^{17}}\right).\label{Landau-1-2}
\end{align}
Combining \eqref{Landau-1} and \eqref{Landau-1-2} gives
\begin{align}
q(n+1)=q_{8}(n+1)+O\left(\frac{1}{n^{16}}\right).
\label{Landau-1-3}
\end{align}
By using the \emph{Mathematica} software, we
expand $q_{8}(n+1)$ into a power series in terms of $n^{-1}$. Noting
formula \eqref{Landau-1-3}, we obtain
\begin{align}
q(n+1)=&q_{8}(n+1)+O\left(\frac{1}{n^{16}}
\label{q8-approximation}\right)\\
=&\frac{1}{n}-\frac 54\frac{1}{n^2}+\frac{49}{32}\frac{1}{n^3}
-\frac{235}{128}\frac{1}{n^4}+\frac{4411}{2048}\frac{1}{n^5}
-\frac{20275}{8192}\frac{1}{n^6}
+\frac{183077}{65536}\frac{1}{n^7}\nonumber\\
&-\frac{815195}{262144}\frac{1}{n^8}+
\frac{28754131}{8388608}\frac{1}{n^9}
-\frac{125799895}{33554432}\frac{1}{n^{10}}
+\frac{1091975567}{268435456}\frac{1}{n^{11}}\nonumber\\
&-\frac{4702048685}{1073741824}\frac{1}{n^{12}}
+\frac{80679143663}{17179869184}\frac{1}{n^{13}}-
\frac{346250976095}{68719476736}\frac{1}{n^{14}}\nonumber\\
&+\frac{2947620308941}{549755813888}\frac{1}{n^{15}}
+O\left(\frac{1}{n^{16}}\right).\nonumber
\end{align}
The above expression is also used in the first and second-correction below. In addition, it is not difficult to obtain
\begin{align}
-\ln(n+\frac 34)+
\ln(n+\frac 74)=\frac{1}{n}-\frac{5}{4}\frac{1}{n^2}
+\frac{79}{48}\frac{1}{n^3}+O\left(\frac{1}{n^{4}}\right).
\label{Landau-1-4}
\end{align}
Inserting \eqref{q8-approximation} and \eqref{Landau-1-4} into
\eqref{L-u0-difference-1} yields
\begin{align}
u_0(n)-u_0(n+1)=\frac{11}{96\pi}\frac{1}{n^3}
+O\left(\frac{1}{n^{4}}\right).\label{Landau-1-5}
\end{align}
Note that the inequalities \eqref{Landau-Mortici} implies $u_0(\infty)=0$. Applying Lemma 1, we obtain
\begin{align}
\lim_{n\rightarrow\infty}n^2u_0(n)=\frac{11}{192\pi}=C_0.
\label{Landau-1-6}
\end{align}

\noindent{\bf (Step 2) The first-correction.} We let
\begin{align}
\eta_1(n)=\frac{C_0}{\Phi_1(2;n)}
=\frac{C_0}{n^2+a_1n+a_0},\label{Landau2-1}
\end{align}
and define
\begin{align}
u_1(n)=G(n)-\eta_0(n)-\eta_1(n).\label{Landau2-2}
\end{align}
Hence
\begin{align}
u_1(n)-u_1(n+1)=&\left(u_0(n)-\eta_1(n)\right)
-\left(u_0(n+1)-\eta_1(n+1)\right)\label{Landau2-3}\\
=&\left(u_0(n)-u_0(n+1)\right)
-\left(\eta_1(n)-\eta_1(n+1)\right).\nonumber
\end{align}
Note that the first term of \eqref{Landau2-3} can be treated
by the same method in (step 1). Here we only need to replace \eqref{Landau-1-4} by the following more accurate power series
expansion
\begin{align}
-\ln(n+\frac 34)+
\ln(n+\frac 74)=&\frac{1}{n}-\frac{5}{4}\frac{1}{n^2}
+\frac{79}{48}\frac{1}{n^3}-\frac{145}{64}\frac{1}{n^4}
+\frac{4141}{1280}\frac{1}{n^5}
-\frac{14615}{3072}\frac{1}{n^6}\label{Landau2-4}\\
&+\frac{205339}{28672}\frac{1}{n^7}
-\frac{179945}{16384}\frac{1}{n^8}
+\frac{10083481}{589824}\frac{1}{n^9}
-\frac{7060405}{262144}\frac{1}{n^{10}}\nonumber\\
&+\frac{494287399}{11534336}\frac{1}{n^{11}}
-\frac{865047235}{12582912}\frac{1}{n^{12}}
+\frac{24221854021}{218103808}\frac{1}{n^{13}}\nonumber\\
&-\frac{84777286235}{469762048}\frac{1}{n^{14}}
+\frac{1186886790259}{4026531840}\frac{1}{n^{15}}
+O\left(\frac{1}{n^{16}}\right).\nonumber
\end{align}
By applying \emph{Mathematica} software again, we have
\begin{align}
&\frac{1}{\Phi_1(2;n)}-\frac{1}{\Phi_1(2;n+1)}\label{Landau2-5}\\
=&\frac{2}{n^3}+\frac{-3-3 a_1}{n^4}+
\frac{4+6a_1+4a_1^2-4 a_0}{n^5}\nonumber \\
&+\frac{-5-10a_1-10a_1^2-5 a_1^3+10a_0+10a_1a_0}{n^6}\nonumber \\
&+\frac{6+15a_1+20a_1^2+15a_1^3+6 a_1^4
-20a_0-30a_1a_0-18a_1^2 a_0+6a_0^2}{n^7}+O\left(\frac{1}{n^{8}}\right).\nonumber
\end{align}
Now combining \eqref{Landau2-3}, \eqref{L-u0-difference-1},\eqref{q8-approximation}, \eqref{Landau2-4}
and \eqref{Landau2-5}, and performing some simplifications, we can obtain
\begin{align}
\pi\left(u_1(n)-u_1(n+1)\right)=&\frac{\frac{235}{128}+\frac{-134 + 11 a_1}{64} }{n^4}\label{Landau2-6}\\
&+\frac{-\frac{4411}{2048}+\frac{11543-1320 a_1-880 a_1^2+880a_0}{3840}}{n^5}\nonumber\\
&+\frac{\frac{20275}{8192}+\frac{
5(-2747+352a_1+352a_1^2+176a_1^3-352a_0-352a_1a_0)}{3072}}{n^6}
\nonumber\\
&+\frac{-\frac{183077}{65536}+\frac{\sigma}{86016}}{n^7}
+O\left(\frac{1}{n^{8}}\right),\nonumber
\end{align}
where
\begin{align*}
\sigma=586449-73920a_1-98560a_1^2-73920a_1^3-29568a_1^4+98560a_0+
147840a_1a_0+88704a_1^2a_0-29568a_0^2.
\end{align*}
The fastest sequence $(u_1(n))_{n\ge 1}$ is obtained when the first
two coefficients of this power series vanish. In this case
\begin{align}
a_1=\frac 32,\quad a_0=\frac{5501}{7040},\label{Landau2-7}
\end{align}
thus
\begin{align*}
u_1(n)-u_1(n+1)=\frac{89684299}{3027763200\pi}\frac{1}{n^7}+
O\left(\frac{1}{n^{8}}\right).
\end{align*}
Finally, by using Lemma 1, one has
\begin{align}
\lim_{n\rightarrow\infty}n^6u_1(n)=
\frac{89684299}{18166579200\pi}=C_1.
\label{Landau2-8}
\end{align}

\noindent{\bf (Step 3) The second-correction.} We let
\begin{align}
\eta_2(n)=\frac{C_1}{\Phi_2(6;n)}
=\frac{C_1}{n^6+b_5n+b_4n^4+b_3n^3+b_2n^2+b_1n+b_0},\label{Landau3-1}
\end{align}
and define
\begin{align}
u_2(n)=G(n)-\eta_0(n)-\eta_1(n)-\eta_2(n).\label{Landau3-2}
\end{align}
Thus
\begin{align}
u_2(n)-u_2(n+1)=\left(u_0(n)-u_0(n+1)\right)
-\left(\eta_1(n)+\eta_2(n)-\eta_1(n+1)-\eta_2(n+1)\right)
.\label{Landau3-3}
\end{align}
We use \eqref{q8-approximation} and \eqref{Landau2-4} to expand $u_0(n)-u_0(n+1)$ into a power
series as in terms of $n^{-1}$. In addition, as mentioned already in Section 3, one can use a similar \emph{Mathematica} program in Example 1 to find
$b_5, b_4, b_3, b_2, b_1$ and $b_0$ in turn. Here we omit the details.
We write
\begin{align}
C_2=-\frac{5691942495934169497683736629269380931519449}{
 65873649616252391923660120676946385934745600\pi}.\label{Landau3-4}
\end{align}
By using Lemma 1 again, it is not very difficult for us to check the following assertion.
\begin{thm} Let $c_0$, $C_2$ be defined by \eqref{c0-definition} and \eqref{Landau3-4} respectively, and
\begin{align}
u_2(n):=G(n)-\left(\frac{1}{\pi}\ln(n+\frac 34)+c_0+\frac{\frac{11}{192\pi}}{\Phi_1(2;n)}
+\frac{\frac{89684299}{18166579200\pi}}
{\Phi_2(6;n)}\right),\label{Landau3-5}
\end{align}
where
\begin{align}
\Phi_1(2;n)=&(n+\frac 34)^2+\frac{1541}{7040},\\
\Phi_2(6;n)=&(n+\frac 34)^6+\frac{1092000370209}{631377464960}(n+\frac 34)^4
-\frac{111862508515629162375}{181198865117870921728}(n+
\frac 34)^2\\
&+\frac{
1824588073050833974528912179250963}{540823069619183303269309779804160}
.\nonumber
\end{align}
Then we have
\begin{align}
&\lim_{n\rightarrow\infty}n^{15}\left(u_2(n)-u_2(n+1)\right)=14C_2,
\label{Landau-3-6-1}\\
&\lim_{n\rightarrow\infty}n^{14}u_2(n)=C_2.
\label{Landau-3-6}
\end{align}
\end{thm}

\begin{rem}
It should be stressed that
that a ``good" initial-correction is very important
for us to accelerate the convergence. In addition, one may study analogous question by choosing different initial-correction.
\end{rem}

The following Theorem tells us how to improve \eqref{Landau-Falaleev} and \eqref{chen-1}.
\begin{thm} Let $c_0$ be defined by \eqref{c0-definition}. Then
for all integer $n\ge 0$, we have
\begin{align}
\frac{C_1}{(n+\frac 32)^6}
<G(n)-\frac{1}{\pi}\ln(n+\frac 34)-c_0-\frac{\frac{11}{192\pi}}{(n+\frac 34)^2+\frac{1541}{7040}}
<\frac{C_1}{(n+\frac 12)^6},\label{theorem4}
\end{align}
where $C_1=\frac{89684299}{18166579200\pi}$.
\end{thm}
\begin{rem}
In fact, Theorem 4 implies that
$u_1(n)$ is a strictly decreasing function of $n$.
\end{rem}
\proof Although the method used in this section is very similar to that in proof of Theorem 1, we would like to give a full proof for the sake of completeness.
First, we can see that the inequalities \eqref{theorem4}
are true for $n=0$. Hence, in the following we only need to prove that these inequalities are also true for $n\ge 1$. To this end, let
\begin{align}
u_1(n)=G(n)-\frac{1}{\pi}\ln(n+\frac 34)-c_0-\frac{\frac{11}{192\pi}}{\Phi_1(2;n)},\label{Landau-2}
\end{align}
it follows easily from \eqref{G-difference-def}
\begin{align}
u_1(n)-u_1(n+1)=&-\frac{1}{\pi}q(n+1)-\frac{1}{\pi}\ln(n+\frac 34)-\frac{\frac{11}{192\pi}}{\Phi_1(2;n)}\label{Landau-3}\\
&+\frac{1}{\pi}\ln(n+\frac 74)
+\frac{\frac{11}{192\pi}}{\Phi_1(2;n+1)}.\nonumber
\end{align}
Let
\begin{align}
f(x)=&-\frac{1}{\pi}q_6(x+1)-\frac{1}{\pi}\ln(x+\frac 34)-\frac{\frac{11}{192\pi}}{\Phi_1(2;x)}
+\frac{1}{\pi}\ln(x+\frac 74)
+\frac{\frac{11}{192\pi}}{\Phi_1(2;x+1)},\label{Landau-4}\\
g(x)=&-\frac{1}{\pi}q_5(x+1)-\frac{1}{\pi}\ln(x+\frac 34)-\frac{\frac{11}{192\pi}}{\Phi_1(2;x)}
+\frac{1}{\pi}\ln(x+\frac 74)
+\frac{\frac{11}{192\pi}}{\Phi_1(2;x+1)}.\label{Landau-5}
\end{align}
From \eqref{Landau-1} and \eqref{Landau-5}, one has
\begin{align}
g(n)<u_1(n)-u_1(n+1)<f(n).\label{Landau-6}
\end{align}
Firstly, we give the lower bound for $g(n)$, and
the upper bound for $f(n)$, respectively.
%in the form of power function
We set $D_1=\frac{89684299}{432537600}$. By using the \emph{Mathematica} software, we easily obtain
\begin{align}
-f'(x)-\frac{D_1}{\pi(x+1)^8}=-\frac{1}{\pi}
\frac{\Psi_1(21;n)}{\Psi_2(30;n)}<0.
\label{Landau-7}
\end{align}
Noting $f(+\infty)=0$ and utilizing \eqref{Landau-7} and Lemma 3, one has
\begin{align}
f(n)&=-\int_{n}^{\infty}f'(x)dx<\int_{n}^{\infty}\frac{D_1}
{\pi(x+1)^8}dx
=\frac{D_1}{7\pi}\frac{1}{(n+1)^7}\label{Landau-8}\\
&<\frac{D_1}{7\pi}\int_{n+\frac 12}^{n+\frac 32}\frac{1}{x^7}dx.\nonumber
\end{align}
Similarly, we can check
\begin{align}
-g'(x)-\frac{D_1}{\pi(x+\frac 32)^8}=\frac{1}{\pi}\frac{\Psi_3(19;n)}{\Psi_4(28;n)}>0.
\label{Landau-9}
\end{align}
Applying $g(+\infty)=0$ and \eqref{Landau-9}, we obtain
\begin{align}
g(n)&=-\int_{n}^{\infty}g'(x)dx>\int_{n}^{\infty}
\frac{D_1}{\pi(x+\frac 32)^8}dx
=\frac{D_1}{7\pi}\frac{1}{(n+\frac 32)^7}\label{Landau-10}\\
&>\frac{D_1}{7\pi}\int_{n+\frac 32}^{n+\frac 52}\frac{1}{x^7}dx.\nonumber
\end{align}
On the other hand, from $u_1(\infty)=0$ and \eqref{Landau-8}, we have
\begin{align}
u_1(n)=&\sum_{m=n}^{\infty}\left(u_1(m)-u_1(m+1)\right)
<\sum_{m=n}^{\infty}\frac{D_1}{7\pi}\int_{m+\frac 12}^{m+\frac 32}\frac{1}{x^7}dx\label{Landau-11}\\
=&\frac{D_1}{7\pi}\int_{n+\frac 12}^{\infty}\frac{1}{x^7}dx
=\frac{D_1}{42\pi}\frac{1}{(n+\frac 12)^6}.\nonumber
\end{align}
Similarly, it follows from \eqref{Landau-10}
\begin{align}
u_1(n)=&\sum_{m=n}^{\infty}\left(u_1(m)-u_1(m+1)\right)
>\sum_{m=n}^{\infty}\frac{D_1}{7\pi}\int_{m+\frac 32}^{m+\frac 52}\frac{1}{x^7}dx\label{Landau-12}\\
=&\frac{D_1}{7\pi}\int_{n+\frac 32}^{\infty}\frac{1}{x^7}dx=\frac{D_1}{42\pi}\frac{1}{(n+\frac 32)^6}.
\nonumber
\end{align}
%Combining \eqref{Landau-11}and \eqref{Landau-12} completes the proof of Theorem 4 for $n\ge 1$. Finally, we check easily that inequalities \eqref{theorem4}
%are true for $n=0$.
This completes the proof of Theorem 4.\qed

\begin{flushleft}

Xiaodong Cao\\
Department of Mathematics and Physics, \\
Beijing Institute of Petro-Chemical Technology,\\
Beijing, 102617, P. R. China \\
 e-mail: caoxiaodong@bipt.edu.cn \\

\bigskip

Hongmin Xu\\
Department of Mathematics and Physics, \\
Beijing Institute of Petro-Chemical Technology,\\
Beijing, 102617, P. R. China \\
e-mail: xuhongmin@bipt.edu.cn

\bigskip
Xu You \\
Department of Mathematics and Physics, \\
Beijing Institute of Petro-Chemical Technology,\\
Beijing, 102617, P. R. China \\
e-mail: youxu@bipt.edu.cn

\end{flushleft}

\end{document}